\documentclass[a4paper,latexcad,twoside,11pt]{article}
\usepackage{amssymb} 
\def\dfrac{\displaystyle\frac}

\begin{document}
\title{\bf Unimodality and genus distributions }
\author{\bf
Liangxia Wan \thanks{\it E-mail  address: lxwan$@$bjtu.edu.cn. {\small Partially supported by NNSFC under Grants
No. $11201024$} }
 \\
  \small\it Department of Mathematics,
Beijing Jiaotong University, Beijing $100044$, China}

\date{}
\maketitle

{\bf Abstract--} New criteria are shown that certain combinations of finite unimodal polynomials are unimodal. 
As applications, unimodality of several polynomial sequences satisfying dependent recurrence relations and their modes are provided. Then unimodality of genus distributions for some ladders and crosses can be determined. As special cases, that of genus distributions for Closed-end ladders, circular ladders, M\"{o}bius ladders and Ringel ladders and their modes are given, which induces the known results for Closed-end ladders.

{\bf Keywords--} Unimodal, log-concave, genus distribution, embedding

\vskip 5mm
 \noindent{\bf $1$. Introduction}
\vskip 5mm

Let $S_i$, $\gamma$ and $\gamma_M$ denote the surface of genus $i$, the minimum genus and maximum genus of a graph $G$ respectively.  Let $g_i(G)$ denote the number of $2$-cell embeddings of $G$ embedded on $S_i$ for $i\ge 0$. The {\it genus distribution }of $G$ is a sequence of numbers followed:
$$g_\gamma(G), g_{\gamma+1}(G),g_{\gamma+2}(G),\cdots,g_{\gamma_M}(G).$$
$f_G(x)=\sum\limits_{i=\gamma}^{\gamma_M}g_i(G)x^i$ is called the {\it genus polynomial } of $G$. In general, it is an NP-complete problem to determine
the genus distribution of a graph, since evaluating $\gamma$ is an NP-complete problem \cite{Th89}.

A finite sequence of nonnegative numbers $\{a_i\}_{i=q}^n$ 
is {\it unimodal} if there exist index $q\le l\le m\le n$ such that $a_q\le a_{q+1}\le\cdots\le a_{l-1}<a_l=\cdots =a_m> a_{m+1}\ge \cdots\ge a_n$ where integers $0\le q\le n$. The corresponding polynomial $\sum\limits_{i=1}^na_ix^i$ is also called {\it unimodal}. $l,l+1,\cdots, m$ are called the {\it modes} of the sequence. If $l=m$, then $m$ is called a {\it peak} of the sequence. The sequence is {\it log-concave} if $a_i^2\ge a_{i-1}a_{i+1}$ for $1\le i\le n-1$. Obviously, if
a sequence is
log-concave, then it is unimodal.
Unimodality problems
arise naturally in many branches of mathematics and have been extensively investigated. See Stanley¡¯s
\cite{St89} and Brenti¡¯s survey articles \cite{Br94} for details. Throughout this paper each sequence of numbers is that of nonnegative numbers. $[a]$ denote the maximum integer not more than $a$ for a nonnegative number $a$.

Unimodality problems for the genus distribution of a graph include unimodality, log-concavity and the reality of zeros of a genus polynomial. Gross, Robbins and Tucker proposed the conjecture that the genus distribution of every graph is log-concave and showed that genus distributions of bouquets are log-concave \cite{GRT89}. Furst, Gross and Statman showed that genus distributions of Closed-end ladders and cobblestone paths are log-concave \cite{FGS89}. Gross, Mansour and Tucker proved genus distributions of some ring-like graphs are log-concave \cite{GMT13}. Stahl proved that all zeros of genus polynomials of tree-like graphs, cobblestone paths, diamond bands and some vertex-forest multijoins are real and negative \cite{St97}.
 Chen, Liu and Wang also considered zeros of genus polynomials of some graphs \cite{LW07,Ch08}. In 2003, Liu conjectured that the genus distribution of each graph is unimodal \cite{Liu03}.
Zhao and Liu verified that the genus distribution of every tree graph is unimodal \cite{ZL06}.

This paper mainly concerns unimodality and log-concavity of polynomial sequences there exist dependent recurrence relations among them. In Section $2$, new criteria are presented to determine unimodality of a sequence of numbers generated from finite unimodal sequences of numbers. In Section $3$, genus distributions for sets of ladder surfaces are shown to be unimodal or even log-concave. Then unimodality of genus distributions for some ladders and crosses can be verified. As special consequences unimodality of genus distributions for Closed-end ladders, circular ladders, M\"{o}bius ladders \cite{Mc87} and Ringel ladders \cite{Te00} are proved and their modes of genus distributions for these graphs are determined in Section $4$, which induces the known results for Closed-end ladders. Section $5$ gives some problems for further study.


\vskip 5mm
\noindent{\bf 2. Criteria }

\vskip 5mm

\noindent{\bf Theorem $2.1$} {\it Suppose that $\{x_i\}_{i=q_1}^{n_1}$ and $\{y_i\}_{i=q_2}^{n_2}$ are unimodal sequences of numbers for $0\le q_1\le n_1$ and $0\le q_2\le n_2$. Let $l_1,l_1+1,\cdots,m_1$ and $l_2,l_2+1,\cdots,m_2$ be the {\it modes }of $\{x_i\}_{i=q_1}^{n_1}$ and $\{y_i\}_{i=q_2}^{n_2}$ for $j=1,2$ and $l_j\le m_j$ respectively, let $r_j$ be non-negative integers and let $a_j>0$. Then the sequence $\{a_1x_{i-r_1}+a_2y_{i-r_2}\}$ is unimodal if and only if it is unimodal for $\min\{l_1+r_1,l_2+r_2\}
\le i\le\max\{m_1+r_1,m_2+r_2\}$. }
\vskip 3mm
\noindent{\bf Proof.} Put $z_i=a_1x_{i-r_1}+a_2y_{i-r_2}$. Clearly, the sequence $\{z_i\}=\{z_i\}_{i=q}^n$ where $q=\min\limits_{1\le j\le 2}\{q_j+r_j\}$ and $n=\max\limits_{1\le j\le 2}\{n_j+r_j\}$. For brevity, it is denoted by  $\{z_i\}$. Put $l=\min\{l_1+r_1,l_2+r_2\}$ and put $m=\max\{m_1+r_1,m_2+r_2\}$. Since $\{x_i\}_{i=p_1}^{n_1}$ and $\{y_i\}_{i=p_2}^{n_2}$ are unimodal and since $a_j>0$ for $j=1,2$, these follow that
$$z_q\le z_{q+1}\le\cdots \le z_{l-1}<z_l=z_m>z_{m+1}\ge\cdots \ge z_{n}.$$
Thus the result is obvious.
 \hskip 4mm $\Box$

\vskip 3mm
If $m-l\le 3$, then it is obvious that $\{a_1x_{i-r_1}+a_2x_{i-r_2}\}$ is unimodal for $l
\le i\le m$ and so the following conclusion holds.
\vskip 3mm
\noindent{\bf Corollary $2.2$}  {\it Suppose that $\{x_i\}_{i=q_1}^{n_1}$ and $\{y_i\}_{i=q_2}^{n_2}$ are unimodal sequences of numbers for $0\le q_1\le n_1$ and $0\le q_2\le n_2$. Let $l_1,l_1+1,\cdots,m_1$ and $l_2,l_2+1,\cdots,m_2$ be the {\it modes }of $\{x_i\}_{i=q_1}^{n_1}$ and $\{y_i\}_{i=q_2}^{n_2}$ for $j=1,2$ and $l_j\le m_j$ respectively, let $r_j$ be non-negative integers and let $a_j>0$. If $\max\limits_{1\le j\le 2}\{m_j+r_j\}-\min\limits_{1\le j\le 2}\{l_j+r_j\}\le 3$, then the sequence $\{a_1x_{i-r_1}+a_2y_{i-r_2}\}$ is unimodal.  }

\vskip 3mm
\noindent{\bf Example $1.$} Suppose that $\{x_i\}_{i=0}^{3}=\{1,3,3,2\}$ and $\{y_i\}_{i=1}^{4}=\{1,2,3,3\}$. Obviously, they are unimodal and $l_1=1,m_1=2$,$l_2=3$ and $m_2=4$. Now consider the sequence $\{x_{i-1}+2y_i\}$. Clearly, $r_1=1$ and $r_2=0$ and then $l=\min\{1+1,3+0\}=2$ and $m=\max\{2+1,4+0\}=4$. Since $m-l=2<3$, $\{x_{i-1}+2y_i\}$ is unimodal. In fact, $\{x_{i-1}+2y_i\}=\{3,7,9,8\}$ is clearly unimodal.
\vskip 3mm
The following result is obtained by using the same technique in the argument of Theorem $2.1$, which is its generalization.

\vskip 3mm
\noindent{\bf Theorem $2.3$} {\it Suppose that $\{x_i^{(1)}\}_{i=q_1}^{n_1}$, $\{x_i^{(2)}\}_{i=q_2}^{n_2}$,$\cdots$,$\{x_i^{(k)}\}_{i=q_k}^{n_k}$ are $k$ unimodal sequences of numbers for $0\le q_j\le n_j$, $1\le j\le k$ and $k\ge 3$. Let $l_j,l_j+1,\cdots,m_j$ be the {\it modes }of $\{x_i^{(j)}\}$ for $1\le j\le k$ and $l_j\le m_j$, let $r_j$ be non-negative integers and let $a_j>0$. Then the sequence $\{\sum\limits_{j=1}^ka_jx_{i-r_j}^{(j)}\}$ is unimodal if and only if it is unimodal for $\min\limits_{1\le j\le k}\{l_j+r_j\}
\le i\le\max\limits_{1\le j\le k}\{m_j+r_j\}$. }
\vskip 3mm

If $m-l\le 3$, then it is obvious that $\sum\limits_{j=1}^ka_jx_{n-r_j}^{(j)}$ is unimodal for $l
\le n\le m$ and so the following conclusion holds.
\vskip 3mm
\noindent{\bf Corollary $2.4$} {\it Suppose that $\{x_i^{(1)}\}_{i=q_1}^{n_1}$, $\{x_i^{(2)}\}_{i=q_2}^{n_2}$,$\cdots$,$\{x_i^{(k)}\}_{i=q_k}^{n_k}$ are $k$ unimodal sequences of numbers for for $0\le q_j\le n_j$, $1\le j\le k$ and $k\ge 3$. Let $l_j,l_j+1,\cdots,m_j$ be the {\it modes }of $\{x_i^{(j)}\}_{i=q_j}^{n_j}$ for $1\le j\le k$ and $l_j\le m_j$, let $r_j$ be non-negative integers and let $a_j>0$. If $\max\limits_{1\le j\le k}\{m_j+r_j\}-\min\limits_{1\le j\le k}\{l_j+r_j\}\le 3$, then the sequence $\{\sum\limits_{j=1}^ka_jx_{i-r_j}^{(j)}\}$ is unimodal. }
\vskip 3mm

Put $r_j=0$ in Theorems $2.1$, $2.3$, Corollaries $2.2$ and $2.4$ and then one obtains the results as follows.

\vskip 3mm
\noindent{\bf Corollary $2.5$} {\it Suppose that $\{x_i^{(1)}\}_{i=q_1}^{n_1}$, $\{x_i^{(2)}\}_{i=q_2}^{n_2}$,$\cdots$,$\{x_i^{(k)}\}_{i=q_k}^{n_k}$ are $k$ unimodal sequences of numbers for for $0\le q_j\le n_j$, $1\le j\le k$ and $k\ge 2$. Let $l_j,l_j+1,\cdots,m_j$ be the {\it modes }of $\{x_i^{(j)}\}$ for $1\le j\le k$ and $l_j\le m_j$ and let $a_j>0$. Then the sequence $\{\sum\limits_{j=1}^ka_jx_i^{(j)}\}$ is unimodal if and only if it is unimodal for $\min\limits_{1\le j\le k}l_j
\le i\le\max\limits_{1\le j\le k}m_j$. }

\vskip 3mm
\noindent{\bf Corollary $2.6$} {\it Suppose that $\{x_i^{(1)}\}_{i=q_1}^{n_1}$, $\{x_i^{(2)}\}_{i=q_2}^{n_2}$,$\cdots$,$\{x_i^{(k)}\}_{i=q_k}^{n_k}$ are $k$ unimodal sequences of numbers for $0\le q_j\le n_j$, $1\le j\le k$ and $k\ge 2$. Let $l_j,l_j+1,\cdots,m_j$ be the {\it modes }of $\{x_i^{(j)}\}$ for $1\le j\le k$ and $l_j\le m_j$ and let $a_j>0$. If $\max\limits_{1\le j\le k}m_j-\min\limits_{1\le j\le k}l_j\le 3$, then the sequence $\{\sum\limits_{j=1}^ka_jx_i^{(j)}\}$ is unimodal. }

\vskip 5mm
\noindent{\bf 3. Unimodality of genus distributions for sets of ladder surfaces  }

\vskip 5mm
Suppose that $a_l$ are distinct letters for $l\ge 1$. Let $R_1^n=a_{k_1}a_{k_2}a_{k_3}\cdots
a_{k_r}$, $R_2^n=a_{k_{r+1}}a_{k_{r+2}}a_{k_{r+3}}$
$\cdots a_{k_n}$,
$R_3^n=a_{t_1}^-a_{t_2}^-a_{t_3}^-\cdots a_{t_s}^-$ and
$R_4^n=a_{t_{s+1}}^-a_{t_{s+2}}^-a_{t_{s+3}}^-\cdots a_{t_n}^-$
where $n\geq k_1>k_2>k_3>\cdots >k_r\geq 1$, $1\leq
k_{r+1}<k_{r+2}<k_{r+3}<\cdots <k_n\leq n$, $n\geq
t_1>t_2>t_3>\cdots >t_s\geq 1$, $1\leq
t_{s+1}<t_{s+2}<t_{s+3}<\cdots
<t_n\leq n$ and $0\leq r,s\leq n$, $k_p\neq k_q$, $t_p\neq t_q$ for $p\neq q$. The sets of ladder surfaces ${\cal S}_j^n$ are given below for
$1\leq j\leq 11$:

\vskip 2mm

 \hskip 3mm ${\cal S}_1^n=\{ R_1^nR_2^nR_3^nR_4^n\}$\hskip 4mm ${\cal S}_2^n=\{
R_1^nR_2^nR_4^nR_3^n\}$ \hskip 4mm ${\cal S}_3^n=\{R_1^nR_3^nR_2^nR_4^n\}$
\vskip 2mm
\hskip 3mm
${\cal S}_4^n=\{aR_1^nR_2^na^-R_3^nR_4^n\}$ \hskip 14mm ${\cal S}_5^n=\{
aR_1^nR_3^na^-R_2^nR_4^n\}$\hskip 4mm
\vskip 2mm
\hskip 3mm
${\cal S}_6^n=\{
aR_1^nR_4^na^-R_2^nR_3^n\}$
\hskip 14mm ${\cal S}_7^n=\{
aR_1^na^-R_3^nR_2^nR_4^n\}$ \hskip 4mm
\vskip 2mm
\hskip 3mm
${\cal S}_8^n=\{ R_1^nR_2^naR_3^na^-bR_4^nb^-\}$ \hskip 8mm ${\cal S}_9^n=\{
R_1^nR_3^naR_2^na^-bR_4^nb^-\}$
\vskip 2mm
\hskip 3mm ${\cal S}_{10}^n=\{
R_1^nR_4^naR_2^na^-bR_3^nb^-\}$ \hskip 7mm ${\cal S}_{11}^n=\{
R_1^naR_2^na^-bR_3^nb^-cR_4^nc^-\}$

\vskip 3mm
Let $g_{i_j}(n)$ denote the number of surfaces of genus $i$ in ${\cal S}_j^n$ for $1\le j\le 11$ and $n\ge 1$. Suppose that $\gamma_j$ and $\gamma_{M_j}$ denote the minimum genus and the maximum genus of surfaces in ${\cal S}_j^n$. Then the {\it genus distribution} of ${\cal S}_j^n$ is the sequence of numbers $\big\{g_{i_j}(n)\big\}_{i=\gamma_j}^{\gamma_{M_j}}$:
$$g_{\gamma_j}(n),g_{{\gamma_j}+1}(n),\cdots, g_{\gamma_{M_j}}(n).$$
Since genus distributions of ${\cal S}_j^n$ are obviously log-concave and their modes are easily to obtained by Theorem $2.5$ of \cite{Wan08} for $1\le j\le 11$ and $1\le n\le 3$, we sometimes omit them.
\vskip 5mm
\noindent{\bf 3.1 Log-concavity of genus distributions of ${\cal S}_j^n$ for $j=1,4,6$}

\vskip 5mm

\noindent{\bf Lemma 3.1}(Theorem $2.5$ of \cite{Wan08}) {\it Suppose that $g_{i_j}(n)$ denotes the number of surfaces of genus $i$ in ${\cal S}_j^n$ for $j=1,6$. Let $C_{n}(i)=\left(
                            \begin{array}{c}
                              n-i\\
                              i
                             \end{array}
                        \right)$.
Then, $$g_{i_1}(n)=
   2^{n+i}\dfrac{2n-3i}{n-i}C_{n}(i) \mbox{ for }0\leq
                        i\leq \Big[\dfrac{n}{2}\Big]\mbox{ and }n\geq
                        1$$}

and $$g_{i_6}(n)=2^{n+i-1}\dfrac{2n-3i+2}{n-i+1}C_{n+1}(i),\mbox{ if }0\leq
                        i\leq \Big[\dfrac{n+1}{2}\Big]\mbox{ and }n\geq 1.$$

\vskip 3mm
\noindent{\bf Theorem 3.2} {\it Let $p_j(n)$ denote peaks of genus distributions of ${\cal S}_j^n$ for $j=1,6$ and $i\ge 1$. Genus distributions of ${\cal S}_j^n$ are log-concave for $j=1,6$ and $n\ge 1$. Then
\begin{equation}
p_1(n)=\Big[\dfrac{n+1}{3}\Big], \mbox{ if }n\ge 1 \mbox{ and }n\ne 2
\end{equation}
\noindent and
\begin{equation}
p_6(n)=\Big[\dfrac{n+2}{3}\Big], \mbox{ if }n\ge 2.
\end{equation}}
\vskip 3mm
\noindent{\bf Proof.} By Lemma 3.1, for $1\le i\le \dfrac{n}{2}$, $$\dfrac{g_{i_1}(n)}{g_{(i-1)_1}(n)}
=2\cdot\dfrac{n+1-2i}{i}\cdot\dfrac{n+2-2i}{n-i}\cdot\dfrac{2n-3i}{2n+3-3i}.$$
Since each item is non-increasing, $\dfrac{g_{i_1}(n)}{g_{(i-1)_1}(n)}$ is non-increasing and therefore the genus distribution of ${\cal S}_1^n$ is log-concave.
%

Now we compute $p_1(n)$. Put $n=3m+k$ for $m\ge 0$ and $k=0,1,2$. We verify the case $k=2$ and leave others to readers. If $k=2$, then  $\Big[\dfrac{n+1}{3}\Big]=m+1.$ By Lemma $3.1$
\begin{eqnarray*}
g_{m_1}(n) & = & 2^{3m+2+m}\dfrac{2(3m+2)-3m}{3m+2-m}C_{3m+2}(m)\\
          &  = & 2^{4m+2}(3m+4)\dfrac{(2m+1)2m(2m-1)\cdots (m+3)}{m!},
          \end{eqnarray*}

\begin{eqnarray*}
g_{(m+1)_1}(n) & = & 2^{3m+2+m+1}\dfrac{2(3m+2)-3(m+1)}{3m+2-(m+1)}C_{3m+2}(m+1)\\
          &  = & 2^{4m+3}(3m+1)\dfrac{2m(2m-1)(2m-2)\cdots (m+2)}{m!}
\end{eqnarray*}
and
\begin{eqnarray*}
g_{(m+2)_1}(n) & = & 2^{3m+2+m+2}\dfrac{2(3m+2)-3(m+2)}{3m+2-(m+2)}C_{3m+2}(m+2)\\
          &  = & 2^{4m+4}(3m-2)\dfrac{(2m-1)(2m-2)\cdots (m-1)}{(m+2)!}.
\end{eqnarray*}

It is clear that
$$
\dfrac{g_{(m+1)_1}(n)}{g_{m_1}} =  \dfrac{6m^2+14m+4}{6m^2+11m+4}> 1$$
and that
$$
\dfrac{g_{(m+1)_1}(n)}{g_{(m+2)_1}}  =  \dfrac{3m^2+7m+2}{3m^2-5m+2}> 1.$$
Thus $\Big[\dfrac{n+1}{3}\Big]$ is the peak of ${\cal S}_1^n$.

For $n\ge 2$, it is easily verified that for $0\leq
                        i\leq \Big[\dfrac{n+1}{2}\Big]$
$$g_{i_6}(n)=\frac{1}{4}g_{i_1}(n+1).$$
Thus the result holds for ${\cal S}_6^n$. \hskip 3mm $\Box$
\vskip 3mm
\noindent{\bf Lemma $3.3$}(Lemma 2.3 of \cite{Wan08})  {\it Let $g_{i_j}(n)$ be the number of
surfaces in $S_n^j$ with genus $i$ for $j=1$ and $4$. Then for
$n\geq 1$ and $1\leq i\leq
\Big[\dfrac{n+1}{2}\Big]$,
$$g_{i_4}(n)=
4g_{(i-1)_1}(n-1).
 $$}

The following result is implied by Theorem $3.2$ and Lemma $3.3$.

\vskip 3mm
\noindent{\bf Corollary $3.4$} {\it Let $p_4(n)$ denote the modes of genus distribution of ${\cal S}_4^n$ for $n\ge 1$. The genus distribution of ${\cal S}_4^n$ is log-concave for each $n$  and

$$
\hskip 17mm p_4(n)=
 \left\{
\begin{array}{ll}
1,2, \mbox{ if }n=3;\\
\Big[\dfrac{n}{3}\Big]+1, \mbox{ otherwise. }
\end{array}
\right.
$$}

\vskip 3mm
\noindent{\bf 3.2 Unimodality of genus distributions of ${\cal S}_j^n$ for $j=2,3,7,8,10$}
\vskip 5mm

\noindent{\bf Lemma 3.5}(Theorem 2.5 of \cite{Wan08}) {\it Let $g_{i_3}(n)$ be the number of
surfaces in $S_n^3$ with genus $i$. Let
$A_{n}(i)=\dfrac{2n-3i-2}{n-2i-1}$,
let $B_{n}(i)=\dfrac{n-i-1}{n-2i}$ and let $C_{n}(i)=\left(
                            \begin{array}{c}
                              n-2-i\\
                              i
                             \end{array}
                        \right).$
Then,
$$
\hskip 17mm g_{i_3}(n)=\left\{
 \begin{array}{llllll}
      2^n+4n-2,\mbox{ if }i=0\mbox{ and }n\geq 1;\\
    C_{n+2}(i+1)\Big(2^{3i+1}A_{n+2}(i+1)\\
       \hskip 15mm+(2^{n+i-1}-2^{3i-2})\dfrac{(i+1)
                        A_{n+2}(i)B_{n+2}(i+1)}{n-2i-1}\Big),\\
                                 \hskip 25mm\mbox{
                        if }1\leq i\leq
                        \Big[\dfrac{n}{2}\Big]-1\mbox{ and }n\geq
                        2;\\
C_{n+1}(i)\Big(2^{3i+1}+(2^{n+i-1}-2^{3i-2})
                        A_{n+2}(i)B_{n+2}(i+1)\Big),\\
                              \hskip 25mm\mbox{
                        if }\Big[\dfrac{n}{2}\Big]-1< i\leq
                        \Big[\dfrac{n-1}{2}\Big]\mbox{ and }n\geq
                        2;\\
(2^{n+i-1}-2^{3i-2})A_{n+2}(i)C_{n+2}(i),\\
                           \hskip 25mm\mbox{
                        if }\Big[\dfrac{n-1}{2}\Big]< i\leq
                        \Big[\dfrac{n}{2}\Big]\mbox{ and }n\geq
                        2.\\
  \end{array}
  \right.
  $$
  }

\noindent{\bf Lemma $3.6$}(Lemma 2.3 of \cite{Wan08}) {\it Let $g_{i_j}(n)$ be the number of
surfaces in $S_n^j$ with genus $i$ for $j=3,7,10$ and
$n\geq 0$. Let $f_{S_j^0}(x)=1$. Then, for
$n\geq 1$,

$$
\hskip 17mm g_{i_j}(n)=
 \left\{
\begin{array}{llllll}
g_{i_3}(n-1)+g_{i_6}(n-1)+2g_{i_7}(n-1),\\
   \hskip 20mm\mbox{ if }j=3,0\leq
i\leq
\Big[\dfrac{n}{2}\Big]\mbox{ and }n\geq 1;\\
2g_{(i-1)_3}(n-1)+2g_{i_{10}}(n-1),\\
     \hskip 20mm\mbox{ if }j=7,0\leq i\leq
\Big[\dfrac{n+1}{2}\Big]\mbox{ and }n\geq 1;\\
g_{(i-1)_6}(n-1)+2g_{(i-1)_7}(n-1)+g_{i_{10}}(n-1),\\
    \hskip 20mm\mbox{ if
}j=10,0\leq i\leq
\Big[\dfrac{n}{2}\Big]+1\mbox{ and }n\geq 1;\\
0, \mbox{ otherwise. }
\end{array}
\right.
 $$}
\vskip 3mm
Lemmas $3.5-6$ imply the following result.

\vskip 3mm

\noindent{\bf Lemma $3.7$} {\it Let $g_{i_j}(n)$ be the number of
surfaces in $S_j^n$ with genus $i$ for $j=3,7,10$ and
$n\geq 1$. Then
$$
g_{i_{7}}(n)=\left\{
\begin{array}{ll}
2, \mbox{ if }i=0;\\
4g_{(i-1)_{3}}(n-1)-2, \mbox{ if }i=1;\\
4g_{(i-1)_{3}}(n-1), \mbox{ otherwise }
\end{array}
\right.
$$
\noindent and
$$
g_{i_{10}}(n)=\left\{
\begin{array}{ll}
1, \mbox{ if }i=0;\\
g_{(i-1)_{3}}(n)-1, \mbox{ if }i=1;\\
g_{(i-1)_{3}}(n), \mbox{ otherwise. }
\end{array}
\right.
$$}

\vskip 3mm
Therefore it is enough to find the unimodality of the genus distribution of ${\cal S}_3^n$ in order to study that of genus distributions of ${\cal S}_7^n$ and ${\cal S}_{10}^n$.
\vskip 3mm

\noindent{\bf Lemma $3.8$} {\it Let $g_{i_3}(n)$ be the number of
surfaces in $S_3^n$ with genus $i$ and
$n\geq 0$. Then for $n\ge 8$,
$$g_{[\frac{n+2}{3}]_3}(n)> g_{([\frac{n+2}{3}]-1)_3}(n)\mbox{ and }
g_{[\frac{n+2}{3}]_3}(n)>g_{([\frac{n+2}{3}]+1)_3}(n). $$}

\vskip 3mm
\noindent{\bf Proof.} Put $n=3m+k$ for $k=0,1,2$. We verify the case $k=0$ and leave others to readers. If $m=3$, then by Lemma 3.5
$$
g_{2_3}(9)=C_{9+2}(2+1)\Big(2^{6+1}A_{9+2}(2+1)+(2^{9+2-1}-2^{6-2})\dfrac{(2+1)
                        A_{9+2}(2)B_{9+2}(2+1)}{9-4-1}\Big)=56432,$$

$$g_{3_3}(9)=C_{9+2}(3+1)\Big(2^{9+1}A_{9+2}(3+1)+(2^{9+3-1}-2^{9-2})\dfrac{(3+1)
                        A_{9+2}(3)B_{9+2}(3+1)}{9-6-1}\Big)=126080$$
\noindent and
$$g_{4_3}(9)=C_{9+1}(4)\Big(2^{13}+(2^{9+4-1}-2^{12-2})
                        A_{9+2}(4)B_{9+2}(4+1)\Big)=69632.$$
Thus $$g_{3_3}(9)>g_{2_3}(9)\mbox{ and }g_{3_3}(9)>g_{4_3}(9).$$
For $m\ge 4$, by Lemma 3.5
\begin{eqnarray*}
g_{m_3}(3m) & = & C_{3m+2}(m+1)\Big(2^{3m+1}A_{3m+2}(m+1)\\
        &  + & (2^{3m+m-1}-2^{3m-2})\dfrac{(m+1)
                        A_{3m+2}(m)B_{3m+2}(m+1)}{3m-2m-1}\Big)\\
                & = & \dfrac{(2m-1)(2m-2)\cdots (m+2)}{(m-1)!}\Big(2^{3m+1}(3m-1)\\
       & + & (2^{4m}-2^{3m-1})(3m+2)\Big)
\end{eqnarray*}
\begin{eqnarray*}
g_{(m-1)_3}(3m) & = & C_{3m+2}(m-1+1)\Big(2^{3(m-1)+1}A_{3m+2}(m-1+1)\\
                & + & (2^{3m+m-1-1}
        -2^{3(m-1)-2})\dfrac{(m-1+1)
                        A_{3m+2}(m-1)B_{3m+2}(m-1+1)}{3m-2(m-1)-1}\Big)\\
        & = & \dfrac{2m(2m-1)\cdots (m+2)}{(m-1)!}\Big(2^{3m-2}\dfrac{3m+2}{m+1}\\
        & + & (2^{4m-2}-2^{3m-5})\dfrac{m(3m+5)(2m+1)}{(m+3)(m+2)(m+1)}\Big)
\end{eqnarray*}
\begin{eqnarray*}
g_{(m+1)_3}(3m) & = & C_{3m+2}(m+1+1)\Big(2^{3(m+1)+1}A_{3m+2}(m+1+1)\\
           & + & (2^{3m+m+1-1}-2^{3(m+1)-2})\dfrac{(m+1+1)
                        A_{3m+2}(m+1)B_{3m+2}(m+1+1)}{3m-2(m+1)-1}\Big)\\
           & = & \dfrac{(2m-2)(2m-3)\cdots (m+3)}{(m-3)!}\Big(2^{3m+4}(3m-4)\\
        & + & (2^{4m}-2^{3m+1})\dfrac{(m+2)(3m-1)(2m-1)}{(m-3)(m-2)(m-1)}\Big)
\end{eqnarray*}
Then
\begin{eqnarray*}
\dfrac{g_{m_3}(3m)}{g_{(m-1)_3}(3m)} & = & \dfrac{2^{3m+1}\dfrac{3m-1}{2m}
        +  (2^{4m}-2^{3m-1})\dfrac{(3m+2)}{2m}}{2^{3m-2}\dfrac{3m+2}{m+1}
         + (2^{4m-2}-2^{3m-5})\dfrac{m(3m+5)(2m+1)}{(m+3)(m+2)(m+1)}}\\
         & \ge & \dfrac{2^{3m+1}
        +  (3\cdot 2^{4m-1}-3\cdot 2^{3m-2})}{3\cdot 2^{3m-2}
         + (3\cdot 2^{4m-1}-3\cdot 2^{3m-4})}\\
         & > & 1
\end{eqnarray*}

\noindent and

\begin{eqnarray*}
\dfrac{g_{m_3}(3m)}{g_{(m+1)_3}(3m)} & = & \dfrac{(2m-1)(m+2)}{(m-2)(m-3)}
\dfrac{2^{3m+1}(3m-1)+(2^{4m}-2^{3m-1})(3m+2)}{2^{3m+4}(3m-4) +  (2^{4m}-2^{3m+1})\dfrac{(m+2)(3m-1)(2m-1)}{(m-3)(m-2)(m-1)}}\\
& > &
\dfrac{2^{3m+2}(3m-1)+(2^{4m+1}-2^{3m})(3m+2)}{ 2^{3m+4}(3m-4) +  (3\cdot 2^{4m+1}-3\cdot 2^{3m+2})}\\
& > & 1.
\end{eqnarray*}
Thus the result holds. \hskip 5mm $\Box$
\vskip 3mm
\noindent{\bf Theorem $3.9$} Let $p_3(n)$ denote the peak of the genus distribution of ${\cal S}_3^n$ for $n\ge 1$. Genus distributions of ${\cal S}_3^n$ are unimodal and
$$
p_3(n)=\left\{
\begin{array}{ll}
\Big[\dfrac{n-1}{2}\Big],\mbox{ if }1\le n\le 6;\\
2,\mbox{ if }n=7;\\
\Big[\dfrac{n+2}{3}\Big],\mbox{ if }n\ge 8.
\end{array}
\right.
$$
\vskip 3mm
\noindent{\bf Proof.} The result is clear by Lemma $3.5$ for $n\ge 8$. Next we verify the result by induction on $n(n\ge 8)$. Assume that the result holds for less than $n(n\ge 9)$.

Now we consider unimodality of the genus distribution of ${\cal S}_3^{n}$. By Lemma $3.6$,
$$g_{i_3}(n)=g_{i_3}(n-1)+g_{i_6}(n-1)+2g_{i_7}(n-1).$$
 Here, ${\cal S}_3^{n-1}$ is unimodal and $p_3(n-1)=\Big[\dfrac{n+1}{3}\Big]$ by the induction hypothesis. ${\cal S}_6^{n-1}$ is unimodal and $p_6(n-1)=\Big[\dfrac{n+1}{3}\Big]$ by Theorem 3.2.
Clearly, ${\cal S}_7^{n-1}$ is unimodal and $p_7(n-1)=\Big[\dfrac{n}{3}\Big]+1$ by Lemma $3.7$ and induction hypothesis.
Then by applying Corollary $2.6$, ${\cal S}_3^{n}$ is unimodal. Combining with Lemma $3.8$, we get $$p_3(n)=\Big[\dfrac{n+2}{3}\Big]$$
as desired. \hskip 3mm $\Box$

Armed with Theorem $3.9$ and Lemma $3.7$, the following conclusion is easily induced.
\vskip 3mm
\noindent{\bf Corollary $3.10$} Let $p_j(n)$ denote the peak of genus distribution of ${\cal S}_j^n$ for $j=7,10$ and $n\ge 1$. Genus distributions of ${\cal S}_j^n$ are unimodal  and
$$
p_7(n)=\left\{
\begin{array}{ll}
\Big[\dfrac{n}{2}\Big],\mbox{ if }2\le n\le 7;\\
3,\mbox{ if }n=8;\\
\Big[\dfrac{n+1}{3}\Big]+1,\mbox{ if }n\ge 9
\end{array}
\right.
$$
\noindent and
$$
p_{10}(n)=\left\{
\begin{array}{ll}
\Big[\dfrac{n+1}{2}\Big],\mbox{ if }1\le n\le 6;\\
3,\mbox{ if }n=7;\\
\Big[\dfrac{n+2}{3}\Big]+1,\mbox{ if }n\ge 8.
\end{array}
\right.
$$
\vskip 3mm

By Lemma $2.3$ of \cite{Wan08} for
$n\geq 1$
$$\hskip 10mmg_{i_j}(n)= \left\{
\begin{array}{ll}
4g_{i_7}(n-1),
   \mbox{ if }j=2,0\leq
i\leq
\Big[\dfrac{n}{2}\Big];\\
4g_{(i-1)_7}(n-1),\mbox{ if }j=8,1\leq i\leq
\Big[\dfrac{n}{2}\Big]+1.
\end{array}
\right.
 $$
Thus
\vskip 3mm
\noindent{\bf Corollary $3.11$} Let $p_j(n)$ denote the peak of the genus distribution of ${\cal S}_j^n$ for $j=2,8$ and $n\ge 3$. Then genus distributions of ${\cal S}_j^n$ are unimodal and
$$
p_2(n)=\left\{
\begin{array}{ll}
\Big[\dfrac{n-1}{2}\Big],\mbox{ if }3\le n\le 8;\\
3,\mbox{ if }n=9;\\
\Big[\dfrac{n}{3}\Big]+1,\mbox{ if }n\ge 10
\end{array}
\right.
$$
\noindent and
$$
p_8(n)=\left\{
\begin{array}{ll}
\Big[\dfrac{n+1}{2}\Big],\mbox{ if }3\le n\le 8;\\
4,\mbox{ if }n=9;\\
\Big[\dfrac{n}{3}\Big]+2,\mbox{ if }n\ge 10.
\end{array}
\right.
$$
\vskip 5mm
\noindent{\bf 3.3 Unimodality of genus distributions of ${\cal S}_j^n$ for $j=5,9,11$}
\vskip 5mm

\noindent{\bf Lemma 3.12}(Theorem 2.5 of \cite{Wan08}) {\it Let $g_{0_5}(1)=2$, $g_{1_5}(1)=2$,
$g_{0_5}(2)=2$, $g_{1_5}(2)=14$,
$g_{0_{9}}(1)=1$,
$g_{1_{9}}(1)=3$, $g_{1_{9}}(2)=10$, $g_{2_{9}}(2)=6$,
$g_{1_{9}}(3)=10$, $g_{2_{9}}(3)=54$,
$B_{n}(i)=\dfrac{n-i-1}{n-2i}$, $C_{n}(i)=\left(
                            \begin{array}{c}
                              n-2-i\\
                              i
                             \end{array}
                        \right)$ and $D_{n}(i)=\dfrac{n}{i}2^{i}$.
Then, $g_{i_j}(n)=$

  $$
   \begin{array}{ll}
                            2^n+8n+8, \mbox{ if }j=5,i=1\mbox{ and }n=3,4;\\
2^n+8n, \mbox{ if }j=5,i=1\mbox{ and }n\geq 5;\\
(2^n-2^{2i-2})C_{n}(i-2)D_n(i-1)+2^{2i}C_{n}(i-1)D_n(i),\\
           \hskip 25mm\mbox{ if }j=5,2\leq i< \dfrac{n}{2}-1\mbox{ and }n\geq 5;\\
(2^n-2^{2i-2})C_{n}(i-2)D_n(i-1)+2^{2i}C_{n}(i-1)D_n(i)+2^{n-1},\\
           \hskip 25mm\mbox{ if }j=5,i= \dfrac{n}{2}-1\mbox{ and }n\geq 5;\\
   \end{array}
$$

$$
\begin{array}{ll}
 (2^n-2^{2i-2})C_{n}(i-2)D_n(i-1)+2^{2i}C_{n}(i-1)D_n(i)+2^n,\\
           \hskip 25mm\mbox{ if }j=5,\dfrac{n}{2}-1< i\leq \dfrac{n-1}{2}\mbox{ and }n\geq 4;\\
 (2^n-2^{2i-2})C_{n}(i-2)D_n(i-1)+2^{\frac{3n}{2}+1}-3\cdot 2^{n-1},\\
           \hskip 25mm\mbox{ if }j=5,\dfrac{n-1}{2}< i\leq \dfrac{n}{2}\mbox{ and }n\geq 4;\\
 (2^n-2^{2i-2})C_{n}(i-2)D_n(i-1),\\
   \hskip 25mm \mbox{ if }j=5,
 \dfrac{n}{2}< i\leq \dfrac{n+1}{2}\mbox{ and }n\geq 3;\\
 6, \mbox{ if }j=9,i=1\mbox{ and }n\geq 4;\\
   3\cdot 2^n+48n-86, \mbox{ if }j=9,i=2\mbox{ and }n=4,5;\\
   3\cdot 2^n+48n-102, \mbox{ if }j=9,i=2\mbox{ and }n\geq 6;\\
  3C_{n}(i-1)\Big(2^{3i-2}B_{n+1}(i)\\
    \hskip 15mm +(2^{n+i-2}-2^{3i-5})\dfrac{(i-1)B_{n}(i-1)B_{n+1}(i-1)}{n-2i+1}\Big),\\
           \hskip 25mm\mbox{ if }j=9,3\leq i< \dfrac{n-1}{2}\mbox{ and }n\geq 6;\\
3C_{n}(i-1)\Big(2^{3i-2}B_{n+1}(i)\\
    \hskip 15mm+(2^{n+i-2}-2^{3i-5})\dfrac{(i-1)B_{n}(i-1)B_{n+1}(i-1)}{n-2i+1}\Big)+2^{n-1},\\
           \hskip 25mm\mbox{ if }j=9,i= \dfrac{n-1}{2}\mbox{ and }n\geq 7;\\
   3C_{n}(i-1)\Big(2^{3i-2}B_{n+1}(i)\\
    \hskip 15mm+(2^{n+i-2}-2^{3i-5})\dfrac{(i-1)B_{n}(i-1)B_{n+1}(i-1)}{n-2i+1}\Big)+2^n,\\
           \hskip 25mm\mbox{ if }j=9,\dfrac{n-1}{2}< i\leq \dfrac{n}{2}\mbox{ and }n\geq 6;\\
 C_{n-1}(i-2)\Big(2^{3i-2}+3(2^{n+i-2}-2^{3i-5})B_{n}(i-1)B_{n+1}(i-1)\Big)\\
    \hskip 15mm+2^{\frac{3n+1}{2}}-3\cdot 2^{n-1},\\
       \hskip 25mm\mbox{ if }j=9,\dfrac{n}{2}<i\leq
       \dfrac{n+1}{2}\mbox{ and }n\geq 5;\\
 3(2^{n+i-2}-2^{3i-5})B_{n+1}(i-1)C_{n}(i-2),\\
    \hskip 25mm \mbox{ if }j=9,
                   \dfrac{n+1}{2}<i\leq \dfrac{n}{2}+1\mbox{ and }n\geq 4.\\
\end{array}
$$
}

\vskip 3mm
Applying a similar way in the argument of Lemma $3.8$, the following conclusion holds.

\vskip 3mm

\noindent{\bf Lemma $3.13$} {\it Let $g_{i_j}(n)$ be the number of
surfaces in $S_j^n$ of genus $i$ for $j=5,9$ and
$n\geq 6$. Then
$$g_{([\frac{n+1}{3}]+1)_5}(n)> g_{[\frac{n+1}{3}]_5}(n)\mbox{ and }
g_{([\frac{n+1}{3}]+1)_5}(n)>g_{([\frac{n+1}{3}]+2)_5}(n)\mbox{ for }6\le n\le 16, $$
$$g_{([\frac{n+2}{3}]+1)_5}(n)> g_{[\frac{n+2}{3}]_5}(n)\mbox{ and }
g_{([\frac{n+2}{3}]+1)_5}(n)>g_{([\frac{n+2}{3}]+2)_5}(n)\mbox{ for } n\ge 17 $$
\noindent and
$$g_{([\frac{n}{3}]+2)_9}(n)> g_{([\frac{n}{3}]+1)_9}(n)\mbox{ and }
g_{([\frac{n}{3}]+2)_9}(n)>g_{([\frac{n}{3}]+3)_9}(n)\mbox{ for }n\ge 10. $$
}

\vskip 3mm

\noindent{\bf Lemma $3.14$}(Lemma 2.3 of \cite{Wan08}) {\it Let $g_{i_j}(n)$ be the number of
surfaces in $S_j^n$ with genus $i$ for $j=5,9,11$ and
$n\geq 0$. Let $f_{S_j^0}(x)=1$. Then, for
$n\geq 1$,
$$\hskip 10mm
g_{i_j}(n)= \left\{
\begin{array}{ll}
2g_{(i-1)_3}(n-1)+2g_{i_9}(n-1),\\
     \hskip 20mm\mbox{ if }j=5,0\leq i\leq
\Big[\dfrac{n+1}{2}\Big]\mbox{ and }n\geq 1;\\
g_{(i-1)_5}(n-1)+2g_{(i-1)_7}(n-1)+g_{i_{11}}(n-1),\\
    \hskip 20mm\mbox{ if
}j=9,1\leq i\leq
\Big[\dfrac{n}{2}\Big]+1\mbox{ and }n\geq 1;\\
2g_{(i-1)_9}(n-1)+2g_{(i-1)_{10}}(n-1),\\
     \hskip 20mm\mbox{ if }j=11,1\leq i\leq
\Big[\dfrac{n+1}{2}\Big]+1\mbox{ and }n\geq 1;\\
0, \mbox{ otherwise. }
\end{array}
\right.
 $$}

\noindent{\bf Theorem 3.15} Let $p_9(n)$ denote the peak of genus distribution of ${\cal S}_9^n$ for $n\ge 2$. Genus distributions of ${\cal S}_9^n$ are unimodal and
$$
p_9(n)=\left\{
\begin{array}{ll}
\Big[\dfrac{n+1}{2}\Big],\mbox{ if }2\le n\le 8;\\
4,\mbox{ if }n=9;\\
\Big[\dfrac{n}{3}\Big]+2,\mbox{ if }n\ge 10
\end{array}
\right.
$$

\vskip 3mm
\noindent{\bf Proof.} We verify the conclusion by induction on $n$. It is obvious for $n=2$. Assume that it holds for less than $n(n\ge 3)$.

 By Lemma 3.14, for $1\leq i\leq
\Big[\dfrac{n}{2}\Big]+1$ and $n\geq 3$,
\begin{equation}g_{i_9}(n)=
g_{(i-1)_5}(n-1)+2g_{(i-1)_7}(n-1)+    2g_{(i-1)_9}(n-2)+2g_{(i-1)_{10}}(n-2).\end{equation}
Here, genus distributions of ${\cal S}_7^{n-1}$ and ${\cal S}_{10}^{n-2}$ are unimodal by Corollary $3.10$. The genus distritution of  ${\cal S}_9^{n-1}$ is unimodal by the induction hypothesis. Now consider unimodality of the genus distribution of ${\cal S}_5^{n-1}$.
By Lemma $3.14$, for $1\leq i\leq
\Big[\dfrac{k}{2}\Big]$ and $k\geq 2$,
$$g_{i_5}(k-1)=
2g_{(i-1)_3}(k-2)+2g_{i_9}(k-2)$$
where the genus distribution of ${\cal S}_3^{n-2}$ is unimodal by Theorem $3.9$ and where
that of  ${\cal S}_9^{n-2}$ is unimodal according to the induction hypothesis. Since it is clear that $p_9(n-2)-(p_3(n-2)+1)\le 3$, that of ${\cal S}_5^{n-1}$ is unimodal by Corollary $2.2$. Combining with Lemma $3.13$, we get for $k\le n$
$$
p_5(k-1)=\left\{
\begin{array}{ll}
\Big[\dfrac{k-1}{2}\Big],\mbox{ if }3\le k\le 6;\\
\Big[\dfrac{k}{3}\Big]+1,\mbox{ if }7\le k\le 17;\\
\Big[\dfrac{k+1}{3}\Big]+1,\mbox{ if }k\ge 18.
\end{array}
\right.
$$
Then we consider the unimodality of genus distribution of ${\cal S}_9^n$.
It is easily known that
\begin{equation}p_7(n-1)\le p_9(n-2)\le p_{10}(n-2)\le p_{5}(n-1)\mbox{ and }p_5(n-1)+1-(p_7(n-1)+1)\le 3.\end{equation}

Thus the conclusion is true for $n$ by armed with ($3-4$), Corollary $2.4$ and Lemma $3.13$ and so is the conclusion by induction.\hskip 5mm $\Box$
\vskip 3mm
By a similar way in the argument of Theorem $3.9$, the following result is obtained.

\vskip 3mm
\noindent{\bf Theorem $3.16$} Let $p_5(n)$ denote the peak of genus distribution of ${\cal S}_5^n$ for each $n\ge 2$. The genus distribution of ${\cal S}_5^n$ is unimodal and
$$
p_5(n)=\left\{
\begin{array}{ll}
\Big[\dfrac{n}{2}\Big],\mbox{ if }2\le n\le 5;\\
\Big[\dfrac{n+1}{3}\Big]+1,\mbox{ if }6\le n\le 16;\\
\Big[\dfrac{n+2}{3}\Big]+1,\mbox{ if }n\ge 17.
\end{array}
\right.
$$

\vskip 3mm
Since Lemma $3.12$ implies for $n\ge 3$
$$
g_{i_{11}}(n)=\left\{
\begin{array}{ll}
2, \mbox{ if }i=1;\\
g_{(i-1)_{5}}(n)-2, \mbox{ if }i=2;\\
g_{(i-1)_{5}}(n), \mbox{ otherwise, }
\end{array}
\right.
$$
it is easy to get the following result.
\vskip 3mm
\noindent{\bf Corollary $3.17$} Let $p_{11}(n)$ denote the peak of the genus distribution of ${\cal S}_j^n$ for each $n\ge 2$. Then the genus distribution of ${\cal S}_{11}^n$ is unimodal and
$$
p_{11}(n)=\left\{
\begin{array}{ll}
\Big[\dfrac{n}{2}\Big]+1,\mbox{ if }2\le n\le 5;\\
\Big[\dfrac{n+1}{3}\Big]+2,\mbox{ if }6\le n\le 16;\\
\Big[\dfrac{n+2}{3}\Big]+2,\mbox{ if }n\ge 17.
\end{array}
\right.
$$
\vskip 3mm

\vskip 5mm
\newpage
 \noindent{\bf $4$. Unimodality of genus distribution of ladders and crosses}
\vskip 5mm

Let $e_0$ and $e_1$ be edges of a connected graph $G_0$. Add vertices $u_1, u_2,\cdots,u_n$ on $e_0$ and add vertices $v_1,v_2,\cdots,v_n$ in sequence for $n\ge 1$.
If one adds $u_lv_l$ denoted by $a_l$ such that they are parallel for $1\le l\le n$, then a ladder $GL_n$ is constructed. Otherwise one adds $u_lv_{n-l+1}$ and then a cross $GC_n$ is obtained(See Fig.1(a) and (b)).

\vskip 8mm


\setlength{\unitlength}{0.97mm}
\begin{center}
\begin{picture}(100,40)
\put(22,32){\circle*{0.6}}
\put(18,32){\circle*{0.6}}
\put(14,32){\circle*{0.6}}
\put(72,36){\circle*{0.6}}
\put(68,36){\circle*{0.6}}
\put(70,36){\circle*{0.6}}
\put(26,40){\circle*{1.5}}
\put(32,40){\circle*{1.5}}
\put(32,24){\circle*{1.5}}
\put(36,40){\line(-1,0){32}}
\put(36,24){\line(-1,0){32}}
\put(86,40){\line(-1,0){32}}
\put(86,24){\line(-1,0){32}}
\put(32,24){\line(0,1){16}}
\put(26,24){\line(0,1){16}}
\put(8,24){\line(0,1){16}}
\put(82,40){\line(-3,-2){24}}
\put(82,24){\line(-3,2){24}}
\put(76,40){\line(-3,-4){12}}
\put(26,24){\circle*{1.5}}
\put(8,24){\circle*{1.5}}
\put(8,40){\circle*{1.5}}
\put(82,40){\circle*{1.5}}
\put(76,40){\circle*{1.5}}
\put(58,40){\circle*{1.5}}
\put(58,24){\circle*{1.5}}
\put(64,24){\circle*{1.5}}
\put(82,24){\circle*{1.5}}
\begin{footnotesize}
\put(32,42){{$u_1$}}
\put(82,42){{$u_1$}}
\put(32,21){{$v_1$}}
\put(82,21){{$v_1$}}
\put(26,42){{$u_2$}}
\put(76,42){{$u_2$}}
\put(8,42){{$u_n$}}
\put(58,42){{$u_n$}}
\put(26,21){{$v_2$}}
\put(8,21){{$v_n$}}
\put(64,21){{$v_{n-1}$}}
\put(58,21){{$v_n$}}
\put(18,16){{$GL_n$}}
\put(70,16){{$GC_n$}}
\put(18,12){{(a)}}
\put(70,12){{(b)}}
\put(30,4){{Fig.1: $GL_n$ and $GC_n$}}
\put(18,42){{$e_0$}}
\put(68,42){{$e_0$}}
\put(18,21){{$e_1$}}
\put(72,21){{$e_1$}}
\end{footnotesize}
\end{picture}
\end{center}

\noindent{\bf Lemma $4.1$} (Theorem $3.1$ of \cite{Wan08})  {\it Let $f_G(x)$ denote the genus polynomial of a graph $G$.
Then
$$f_{GL_n}(x)=\sum\limits_{j=1}^{11}f_j(x)f_{{\cal S}_j^n}(x)$$
\noindent where $f_{G}(x)=\sum\limits_{i=1}^{11}f_j(x)$.}

\vskip 3mm

Similarly, the result is clear for a cross $GC_n$.

\vskip 3mm

\noindent{\bf Lemma $4.2$} {\it Let $f_G(x)$ denote the genus polynomial of a graph $G$. Then
$$f_{GC_n}(x)=\sum\limits_{j=1}^{11}f_j(x)f_{{\cal S}_j^n}(x)$$
\noindent where $f_{G}(x)=\sum\limits_{i=1}^{11}f_j(x)$.}

\vskip 3mm

Based on Lemmas $4.1-2$, unimodality of genus distributions of some ladders and crosses can be determined by using criteria in Section $2$.

\vskip 3mm

\noindent{\bf Lemma $4.3$} (Propositions $2-3$ of \cite{St89}) {\it Let $a(x)$ and $b(x)$ be polynomials with positive coefficients.

 $(1)$ If both $a(x)$ and $b(x)$ are log-concave, then so is $a(x)b(x)$.

 $(2)$ If both $a(x)$ is log-concave and $b(x)$ is unimodal, then $a(x)b(x)$ is unimodal.}
\vskip 3mm
Armed with Lemma $4.3$, unimodality of genus distributions of some ladders and crosses is determined as follows.

\vskip 3mm

\noindent{\bf Theorem $4.4$} {\it Let $G_n$ be a ladder or a cross.
Suppose that $f_{{\cal P}_{G_n}}(x)=f_{{\cal P}_{G_0}}(x)f_{{\cal S}_j^n}(x)$ for $j=1$ or $6$. If $f_{{\cal P}_{G_0}}(x)$ is log-concave$($or unimodal$)$, then $f_{{\cal P}_{G_n}}(x)$ is log-concave $($or unimodal$)$.}

\vskip 3mm

Next we consider several types of ladders which are Closed-end ladders $L_n$, circular ladders $CL_n$, M\"{o}bius ladders $ML_n$, Ringel ladders $RL_n$ and a type of crosses $R_n$. See Fig.$2$ for $n=4$.

\bigskip
\setlength{\unitlength}{0.97mm}
\begin{center}
\begin{picture}(184,58)
\put(33,38){\oval(32,28)}
\put(24,52){\circle*{1.5}}
\put(30,52){\circle*{1.5}}
\put(36,52){\circle*{1.5}}
\put(42,52){\circle*{1.5}}
\put(24,24){\circle*{1.5}}
\put(30,24){\circle*{1.5}}
\put(36,24){\circle*{1.5}}
\put(42,24){\circle*{1.5}}
\put(24,52){\line(0,-1){27.5}}
\put(30,52){\line(0,-1){27.5}}
\put(36,52){\line(0,-1){27.5}}
\put(42,52){\line(0,-1){27.5}}
\put(88,38){\oval(28,28)}
\put(88,38){\circle{12}}
\put(88,24){\circle*{1.5}}
\put(88,32.5){\circle*{1.5}}
\put(88,44){\circle*{1.5}}
\put(88,52){\circle*{1.5}}
\put(74,38){\circle*{1.5}}
\put(82.3,38){\circle*{1.5}}
\put(93.7,38){\circle*{1.5}}
\put(102,38){\circle*{1.5}}
\put(88,24){\line(0,1){8}}
\put(88,52){\line(0,-1){8}}
\put(74,38){\line(1,0){9}}
\put(102,38){\line(-1,0){9}}
\put(130,38){\circle*{1.5}}
\put(122,38){\circle*{1.5}}
\put(136,44){\circle*{1.5}}
\put(136,52){\circle*{1.5}}
\put(142,38){\circle*{1.5}}
\put(150,38){\circle*{1.5}}
\put(136,24){\circle*{1.5}}
\put(136,32){\circle*{1.5}}
\qbezier(142,38)(138,24)(136,24)
\qbezier(150,38)(145,30)(136,32)
\put(130,38){\line(1,-1){6}}
\put(122,38){\line(1,-1){14}}
\put(136,52){\line(1,-1){14}}
\put(136,44){\line(1,-1){6}}
\put(122,38){\line(1,1){14}}
\put(130,38){\line(1,1){6}}
\put(60,0){\oval(28,12)}
\put(52,6){\circle*{1.5}}
\put(57.5,6){\circle*{1.5}}
\put(63,6){\circle*{1.5}}
\put(68,6){\circle*{1.5}}
\put(52,-6){\circle*{1.5}}
\put(57.5,-6){\circle*{1.5}}
\put(63,-6){\circle*{1.5}}
\put(68,-6){\circle*{1.5}}
\put(52,6){\line(0,-1){12}}
\put(57.5,6){\line(0,-1){12}}
\put(63,6){\line(0,-1){12}}
\put(68,6){\line(0,-1){12}}
\put(46,0){\circle*{1.5}}
\put(74,0){\circle*{1.5}}
\qbezier(45.3,0)(60,35)(74.7,0)
\put(120,4){\oval(32,22)}
\put(112,15){\circle*{1.5}}
\put(117.5,15){\circle*{1.5}}
\put(123,15){\circle*{1.5}}
\put(128,15){\circle*{1.5}}
\put(112,-7){\circle*{1.5}}
\put(117.5,-7){\circle*{1.5}}
\put(123,-7){\circle*{1.5}}
\put(128,-7){\circle*{1.5}}
\put(112,15){\line(3,-4){16}}
\put(117.5,15){\line(1,-4){5.5}}
\put(123,15){\line(-1,-4){5.5}}
\put(128,15){\line(-3,-4){16}}
\put(107,13.5){\circle*{1.5}}
\put(107,-5.5){\circle*{1.5}}
\put(107,13.5){\line(0,-1){19}}
\begin{footnotesize}
\put(32,20){{$L_4$}}
\put(85,20){{$CL_4$}}
\put(132.0,20.0){$ML_4$}
\put(57.0,-10){$RL_4$}
\put(118.0,-11){$R_4$}
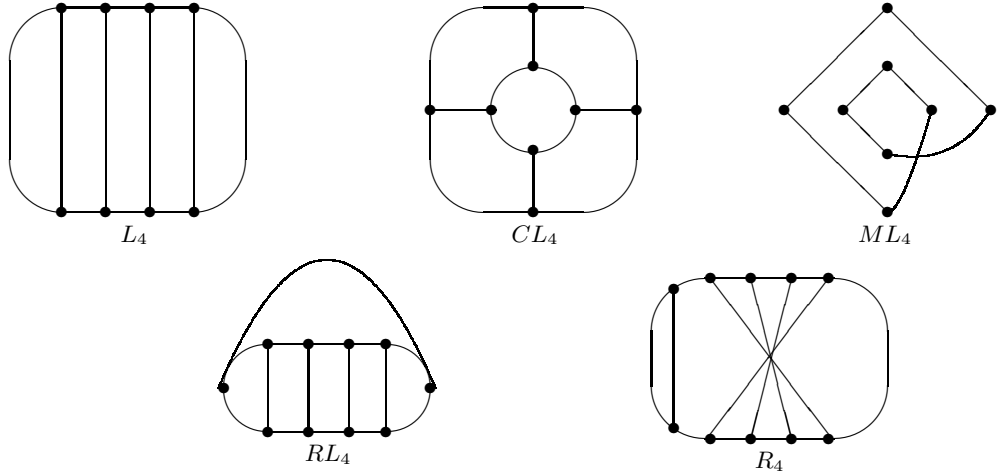
\put(68.0,-15){Fig.$2$:$L_4$,$CL_4$,$ML_4$ and $R_4$}
\end{footnotesize}
\end{picture}
\end{center}
\vskip 16mm

By applying Section $3$ of \cite{Wan08}, $$f_{L_n}(x)=f_{S_6^n}(x).$$ Then the following conclusion is immediate by Theorem $3.2$, which induces the known results for Closed-end ladders in \cite{FGS89}.

\vskip 3mm

\noindent{\bf Corollary $4.5$} {\it Let $p_{L}(n)$ denote the peak of the genus distribution of a Closed-end ladder $L_n$ for each $n\ge 1$. The genus distribution of $L_n$ is log-concave and
$$p_L(n)=\Big[\dfrac{n+2}{3}\Big], \mbox{ if }n\ge 2.$$}


For $CL_n$ and $ML_n$, since genus distribution for $ML_n$ equals to that of $CL_n$, except that $ML_n$ has
two extra embeddings of genus $1$ and two fewer embeddings of
genus $0$(\cite{Mc87}), it is enough to consider unimodality of the genus distribution for $CL_n$. In \cite{Wan08} we have $$g_i(CL_n)=2g_{i_9}(n-1)+2g_{i_{10}}(n-1).$$ Since by Lemma $3.14$
$$g_{i_{11}}(n)=2g_{(i-1)_9}(n-1)+2g_{(i-1)_{10}}(n-1),$$
$$g_{i}(CL_n)=g_{(i+1)_{11}}(n).$$
Thus, we have the following result by Corollary $3.17$.
\vskip 3mm

\noindent{\bf Corollary $4.6$} {\it Let $p_{CL}(n)$ and $p_{ML}(n)$ denote the peaks of genus distributions of circular ladders $CL_n$ and M$\ddot{o}$bius ladders $ML_n$ for $n\ge 2$,respectively. Then their genus distributions are unimodal and
$$p_{CL}(n)=p_{ML}(n)=\left\{
\begin{array}{ll}
\Big[\dfrac{n}{2}\Big],\mbox{ if }2\le n\le 5;\\
\Big[\dfrac{n+1}{3}\Big]+1,\mbox{ if }6\le n\le 16;\\
\Big[\dfrac{n+2}{3}\Big]+1,\mbox{ if }n\ge 17.
\end{array}
\right.
$$}

Similarly,
since  in \cite{Wan08} $$g_i(RL_n)=2g_{(i-1)_3}(n)
+2g_{i_{10}}(n)$$ and since by Lemma 3.6
$$g_{i_7}(n+1)=2g_{(i-1)_3}(n)
+2g_{i_{10}}(n),$$ the following result is implied by Corollary $3.10$.

\vskip 3mm
\noindent{\bf Corollary $4.7$} {\it Let $p_{RL}(n)$ denote the peak of the genus distribution of a Ringel ladder $RL_n$ for each $n\ge 1$. Then its genus distribution is unimodal and
$$
p_{RL}(n)=\left\{
\begin{array}{ll}
\Big[\dfrac{n+1}{2}\Big],\mbox{ if }1\le n\le 6;\\
3,\mbox{ if }n=7;\\
\Big[\dfrac{n+1}{3}\Big]+1,\mbox{ if }n\ge 8
\end{array}
\right.
$$

}

Now consider the type of crosses $R_n$. The following equation is immediate from \cite{WL06} and Theorem $3$ in \cite{WL09}.

$$g_{i}(R_n)=2\mu_{i_6}(n)+2\mu_{{(i-1)}_1}(n)=2g_{i_5}(n)+2g_{{(i-1)}_2}(n).$$
By applying the same technique in the argument of Theorem $3.9$, the following result is obtained.
\vskip 3mm
\noindent{\bf Theorem $4.8$} {\it Let $p_{R}(n)$ denote the peak of the genus distribution of $R_n$ for each $n\ge 1$. Then its genus distribution is unimodal and
$$
p_{R}(n)=\left\{
\begin{array}{ll}
\Big[\dfrac{n+1}{2}\Big],\mbox{ if }1\le n\le 6;\\
3,\mbox{ if }n=7;\\
\Big[\dfrac{n+2}{3}\Big]+1,\mbox{ if }n\ge 8.
\end{array}
\right.
$$

}
\newpage
\vskip 5mm
\noindent{\bf $5$. Further study  }

\vskip 5mm

\noindent{\bf Problem $5.1$} {\it Determine whether genus distributions of sets of ladder surfaces ${\cal S}_j^n$ are log-concave for $2\le j\le 11$, $j\ne 4,6$ and $n\ge 2.$   }
\vskip 3mm

\noindent{\bf Problem $5.2$} {\it Let $k$ sequences of polynomials $\{P_j(n)\}$ satisfy certain dependent recurrence relations for $k\ge 2$ and $1\le j\le k$. If their explicit expressions are unknown, then determine whether they are unimodal (or log-concave). For example, let $P_1(n)=\sum\limits_{i=0}^n g_i(n)x^i$ and let $P_2(n)=\sum\limits_{i=1}^n g_{i_2}(n)x^i$ where $g_{0}(0)=1$, $g_{0}(1)=2$,
$g_{1}(1)=14$, $g_{1_2}(1)=4$ (Theorem $4.1$ \cite{WL08}). For $n\geq 2$,
$$
   \left\{
     \begin{array}{llll}
        g_i(n)= 2g_i({n-1})+8g_{i-1}({n-1})+48g_{i-1}({n-2})+12g_{{(i-1)}_2}({n-1}),\\
                g_{i_2}(n)=8g_{{(i-1)}_2}({n-1})+32g_{i-1}({n-2}).
     \end{array}
   \right.
$$
   }
determine whether $\{P_1(n)\}$ and $\{P_2(n)\}$ are unimodal ( or log-concave ).
\vskip 3mm
\noindent{\bf Problem $5.3$} {\it Suppose that $\{x_i\}_{i=q_1}^{n_1}$ and $\{y_i\}_{i=q_2}^{n_2}$ are unimodal sequences of numbers for $0\le q_1\le n_1$ and $0\le q_2\le n_2$. Let $l_1,l_1+1,\cdots,m_1$ and $l_2,l_2+1,\cdots,m_2$ be the {\it modes }of $\{x_i\}_{i=q_1}^{n_1}$ and $\{y_i\}_{i=q_2}^{n_2}$ for $j=1,2$ and $l_j\le m_j$ respectively, let $r_j$ be non-negative integers and let $a_j>0$. If $\max\limits_{1\le j\le 2}\{m_j+r_j\}-\min\limits_{1\le j\le 2}\{l_j+r_j\}\ge 4$, then determine the conditions such that the sequence of numbers $\{a_1x_{i-r_1}+a_2y_{i-r_2}\}$ is unimodal.  }

\vskip 3mm
\noindent{\bf Problem $5.4$} {\it Suppose that $\{x_i^{(1)}\}_{i=q_1}^{n_1}$, $\{x_i^{(2)}\}_{i=q_2}^{n_2}$,$\cdots$,$\{x_i^{(k)}\}_{i=q_k}^{n_k}$ are $k$ unimodal sequences of numbers for $k\ge 3$. Let $l_j,l_j+1,\cdots,m_j$ be the {\it modes }of $\{x_i^{(j)}\}$ for $1\le j\le k$ and $l_j\le m_j$, let $r_j$ be non-negative integers and let $a_j>0$. If $\max\limits_{1\le j\le 2}\{m_j+r_j\}-\min\limits_{1\le j\le 2}\{l_j+r_j\}\ge 4$, then determine the conditions such that the sequence of numbers $\{\sum\limits_{j=1}^ka_jx_{i-r_j}^{(j)}\}$ is unimodal.  }

\end{document}